\documentclass{article}
\usepackage{tikz}
\usepackage{latexsym,amsfonts,amssymb}
\usepackage[latin1]{inputenc}
\usepackage{amsmath}
\usepackage{mathtools}

\title{\sc Triangles in the graph of conjugacy classes of normal subgroups}

\author{Antonio Beltr\'an\\
\footnotesize
Departamento de Matem\'aticas,\\
\footnotesize Universidad Jaume I, \footnotesize
12071 Castell\'on, Spain\\
\footnotesize
e-mail: abeltran@mat.uji.es\\
\\
Mar\'{\i}a Jos\'e Felipe\\
\footnotesize
Instituto Universitario de Matem\'atica Pura y Aplicada,\\
\footnotesize Universidad Polit\'ecnica de Valencia, \footnotesize
46022 Valencia, Spain\\
\footnotesize e-mail: mfelipe@mat.upv.es \\
\\ Carmen Melchor\\
\footnotesize
Departamento de Educaci\'on,\\
\footnotesize Universidad Jaume I, \footnotesize
12071 Castell\'on, Spain\\
\footnotesize
e-mail: cmelchor@uji.es      }

\date{}
\begin{document} \maketitle

\begin{abstract}
Let $G$ be a finite group and $N$ a normal subgroup of $G$. We determine the structure of $N$ when the graph $\Gamma_{G}(N)$, which is the graph associated to the conjugacy classes of $G$ contained in
$N$, has no triangles and when the graph consists in exactly one triangle.



\end{abstract}

\bigskip

\section{Introduction}
Let $G$ be a finite group and let $N$ be a normal subgroup of $G$ and let $x\in N$. We denote by $x^{G}=\lbrace x^{g} \mid g \in G \rbrace$ the $G$-conjugacy class of $x$. Let $\Gamma_{G}(N)$ be the
graph associated to these $G$-conjugacy classes, which was defined in \cite{Nuestro1} as follows: its vertices are the $G$-conjugacy classes of $N$ of cardinality bigger than 1, that is, $G$-classes
of elements in $N\setminus ({\bf Z}(G)\cap N)$, and two of them are joined by an edge if their sizes are not coprime.\\

The above graph extends the ordinary graph, $\Gamma(G)$, which was formerly defined in \cite{BerHerMann}, whose vertices are the non-central conjugacy classes of $G$ and two vertices are joined by an
edge if their sizes are not coprime. The graph $\Gamma_G(N)$ can be viewed as the subgraph of $\Gamma(G)$ induced by those vertices of $\Gamma(G)$ which are contained in $N$.
This fact does not allow us to obtain directly properties of the graph of $G$-classes.

On the other hand, the graph $\Gamma_G(N)$ is not related to $\Gamma(N)$. In fact, the action of a finite group $G$ on a normal subgroup $N$ determines that the conjugacy classes  of $N$ are grouped
into $G$-classes  and it is possible to find primes dividing some $G$-class sizes of $N$ which do not divide the order of $N$. However, the study of these classes and their sizes has emerged as a new
useful tool to obtain information regarding normal subgroups from a new point of view. 

In this paper we analyze the structure properties of $N$ in two particular cases: when $\Gamma_{G}(N)$ has no triangles and when it consists in exactly one triangle. In the following theorem,
 we extend known results appeared  for the graph $\Gamma(G)$ when it has no triangles. In fact, in \cite{FangZhang} Fang and Zhang gave a complete list of those groups $G$ whose graph $\Gamma(G)$  has
 no triangles. These groups are exactly: the symmetric group $S_{3}$, the dihedral group $D_{10}$,  the three pairwise non-isomorphic non-abelian groups of order 12, and the non-abelian group of order
 21.   In order to obtain our result we firstly need to study the structure of $N$ when $\Gamma_{G}(N)$ has few vertices: exactly one, two or three. The structure of normal subgroups which are union of three or four $G$-conjugacy classes already appeared in
\cite{Sha} and \cite{Riese}, respectively. However, in these papers central $G$-classes are also contemplated, and we need to exclude them.  We remark that some of the results that we obtain
 for $\Gamma_{G}(N)$ are not possible for $\Gamma(G)$.   \\

{\bf Theorem A}. {\it Let $N$ be a non-central normal subgroup of a finite group $G$ such that $\Gamma_{G}(N)$ has no triangles. Then $N$ is a $\{p,q\}$-group and satisfies one of these properties
\begin{enumerate}
\item $N$ is a $p$-group.
\item $N = P \times Q$ with $P$ a $p$-group and  $Q \subseteq {\bf{Z}}(G) \cap N$, $Q \cong { \mathbb Z}_2$.
\item $N=P \times Q$  with $P$ a $p$-group and $Q$ a $q$-group both elementary abelian with $p$ and $q$ odd primes. In this case ${\bf{Z}}(G) \cap N=1$.
\item $N$ is a quasi-Frobenius group with abelian kernel and complement and ${\bf{Z}}(G) \cap N \cong {\mathbb Z}_2$.
\item $N$ is a Frobenius group with complement isomorphic to ${\mathbb Z}_q$, ${\mathbb Z}_{q^2}$ or $Q_8$. In the first case, the kernel of $N$ is a $p$-group with exponent less or equal than
    $p^2$ and in the two latter cases, the kernel of $N$ is $p$-elementary abelian.

\end{enumerate}
}

\bigskip
We recall that a group $G$ is said to be quasi-Frobenius if $G/{\bf Z}(G)$ is Frobenius and then the preimage in $G$ of its kernel and complement are called the kernel and complement of $G$.

\bigskip
Returning to the ordinary case, the structure of a group $G$ having exactly three non-central classes forming a triangle in $\Gamma(G)$ is very limited. We give a proof here.

\bigskip

{\bf Theorem B.} {\it  Let $G$ be a finite group such that $\Gamma(G)$ consists exactly of one triangle. Then $G\cong Q_8$ or $G\cong D_8$.}

\bigskip
We study the structure of a normal subgroup $N$ when $\Gamma_G(N)$  consists in exactly one triangle.

\bigskip
{\bf Theorem C.} {\it Let $N$ be a normal subgroup of a finite group $G$. If $\Gamma_{G}(N)$ has exactly one triangle, then  one of the following possibilities holds:
\begin{enumerate}
\item $N$ is a $p$-group for some prime $p$.
\item $N=P\times Q$, with $P$ $p$-elementary abelian and $Q$ $q$-elementary abelian for some primes $p$ and $q$,  and ${\bf Z}(G)\cap N=1$.
\item $N=P\times Q$, with $P$ a $p$-group for a prime $p\neq 3$,  and $Q\subseteq  {\bf Z}(G)\cap N$, $Q \cong {\mathbb Z}_{3}$ and $P/({\bf Z}(G)\cap P)$ has exponent $p$.
\item $N=PQ$, where $P$ is a Sylow $p$-subgroup, $p\neq 2$ and $Q$ is a Sylow $2$-subgroup of $N$. In addition, $P$ has exponent $p$, $|${\rm \textbf{Z}}$(G)\cap N|=2$ and $Q/(${\rm
    \textbf{Z}}$(G)\cap N)$ is $2$-elementary abelian.
\item Either $N$ is a Frobenius group with complement ${\mathbb Z}_{q}$, ${\mathbb Z}_{q^{2}}$ or $Q_{8}$ for a prime $q$, or there are two primes $p$ and $q$ such that $N/${\rm
    \textbf{O}}$_{p}(N)$ is a Frobenius group of order $pq$ and {\rm \textbf{O}}$_{p}(N)$ has exponent $p$. In this case, ${\bf Z}(G)\cap N =1$.
\item $N\cong A_{5}$ and $G=(N\times K) \langle x \rangle$ for some  $K\leq G$ and $x\in G$, with $x^{2}\in N\times K$ and $G/K\cong N\langle x \rangle \cong S_{5}$.
\end{enumerate}}

We will provide examples for all the cases with the help of \cite{gap}. In particular, we  see that the still open $S_3$-conjecture, which asserts that any finite group in which distinct conjugacy
classes have
distinct sizes is isomorphic to the symmetric group $S_3$, does not hold for $G$-conjugacy classes lying in a normal subgroup. All groups are supposed to be finite.

\section{Preliminary results}
In this section we state some preliminary results which are necessary for the proofs. In \cite{Nuestro1}, the authors obtained the following three theorems, which concern the case where the graph
$\Gamma_{G}(N)$ is disconnected. We denote by $n(\Gamma_{G}(N))$ the number of connected components of $\Gamma_{G}(N)$. \\

{\bf Theorem 2.1.} {\it Let $G$ be group and let $N$ be a normal subgroup of $G$, then $n(\Gamma_{G}(N))\leq 2$.}\\

{\bf Theorem 2.2.} {\it Let $G$ be a group and let $N$ be a normal subgroup of $G$. If $\Gamma_{G}(N)$ has two connected components then, either $N$ is quasi-Frobenius with abelian kernel and
complement, or $N=P\times A$ where $P$ is a $p$-group and $A\leqslant$ {\rm \textbf{Z}}$(G)$.}\\

{\bf Theorem 2.3.} {\it Let $G$ be a group and let $N$ be a normal subgroup of $G$. If $n(\Gamma_{G}(N))=2$, each of its connected components is a complete graph.}\\

It will be hugely important to use the following two results concerning CP-groups, that is, those  groups having all elements of prime-power order, due to Higman (\cite{Higman}) and Deaconescu
(\cite{Deaconescu}), respectively. \\

{\bf Theorem 2.4.} {\it Let $G$ be a solvable group all of whose elements have prime power order. Let $p$ be the prime such that $G$ has a normal $p$-subgroup greater than 1, and let $P$ be the
greatest normal $p$-subgroup of $G$. Then  $G=P$ or $G/P$ is either
\begin{enumerate}
\item a cyclic group whose order is a power of a prime other than $p$; or
\item a generalized quaternion group, $p$ being odd; or
\item a group of order $p^{a}q^{b}$ with cyclic Sylow subgroups, $q$ being a prime of the form $kp^{a}+1$.
\end{enumerate} Thus, $G$ has order divisible by at most two primes, and $G/P$ is metabelian.}\\

We note that in cases $1$ and $2$ of the above theorem, $G$ is a Frobenius group since it is immediate that a $p$-complement of $G$ acts fixed-point-freely on $P$. In  $3$, however, $G/P$ is a
Frobenius group. \\

{\bf Theorem 2.5.} {\it Let $G$ be a group having all nontrivial elements of prime order. Then one of the following cases occurs:
\begin{enumerate}
\item $G$ is a $p$-group of exponent $p$.
\item \begin{enumerate}
\item $|G|=p^{a}q$, $3\leq p < q$, $a\geq 3$, $|F(G)|=p^{a-1}$, $|G:G'|=p.$
\item $|G|=p^{a}q$, $3\leq q < p$, $a\geq 1$, $|F(G)|=|G'|=p^{a}.$
\item $|G|=2^{a}p$, $p\geq 3$, $a\geq 2$, $|F(G)|=|G'|=2^{a}.$
\item $|G|=2p^{a}$, $p\geq 3$, $a\geq 1$, $|F(G)|=|G'|=p^{a}.$ and $F(G)$ is elementary abelian.
\end{enumerate}
\item $G\cong A_{5}.$
\end{enumerate}}

\section{$\Gamma_{G}(N)$ with one vertex}

We start with some properties of those normal subgroups having exactly one non-central $G$-conjugacy class. We stress, however, that for any group $G$, it is elementary that the case in which
$\Gamma(G)$ has exactly one vertex cannot happen. \\

{\bf Theorem 3.1.} {\it If $G$ is a group and $N$ is a normal subgroup of $G$ such that $\Gamma_{G}(N)$ has only one vertex, then $N$ is a $p$-group for some prime $p$ and $N/(N\cap{\bf Z}(G))$ is an
elementary abelian $p$-group.}

\bigskip

{\it Proof.}  Since $\Gamma_{G}(N)$ has one vertex, the set $N\setminus (N\cap \textbf{Z}(G))$  consists of exactly one $G$-conjugacy class and so, $N/(N\cap{\rm \textbf{Z}}(G))$ is characteristically
simple. Thus, $N/(N\cap \textbf{Z}(G))$ is $p$-elementary abelian for some prime $p$. Moreover, $N \cap \textbf{Z}(G)$ is also a $p$-group, and as a consequence, $N$ is a $p$-group too. In fact, if a
prime $q\neq p $ divides $|N \cap \textbf{Z}(G)|$, then $N/(N\cap{\rm \textbf{Z}}(G))$ would have elements of order divisible by $p$ and $q$, a contradiction. $\Box$

\bigskip

{\bf Example 1.} Let $G=$ SL(2,3) and let $N$ be the normal subgroup of $G$ isomorphic to $Q_{8}$. Then $\Gamma_{G}(N)$ has exactly one vertex of size 6 and $N/({\rm \textbf{Z}}(G)\cap N)$ has
exponent 2.

\section{$\Gamma_{G}(N)$ with two vertices}

We divide this section into two cases: when there are exactly two vertices joined by an edge, and when they are not joined.
\subsection{$\Gamma_{G}(N)$ with two vertices without edge}
It was proved in \cite{BerHerMann} that if $\Gamma(G)$ has vertices but no edges, then $G$ is isomorphic to $S_{3}$, and this is clearly not true for $\Gamma_{G}(N)$, just take for instance $G=S_3$
and $N=A_3$ or Example 1.  The following property helps us to further refine the structure of these normal subgroups.\\

{\bf Lemma 4.1.1.} {\it Let $N$ be a normal subgroup of a group $G$. If $\Gamma_{G}(N)$ has two vertices but no edges, then {\rm \textbf{Z}}$(G)\cap N=1$.}\\

{\it Proof.} Let $A$ and $B$ be the two vertices of $\Gamma_{G}(N)$ and let $z \in \textbf{Z}(G)\cap N$. Then $zA$ is a $G$-class of $N$ of cardinality $|A|$. Hence, $zA=A$ and $\langle z \rangle
A=A$. Thus, $|\langle z \rangle|$ divides $|A|$. Similarly, we deduce that $|\langle z\rangle |$ divides $|B|$ and since $(|A|,|B|)=1$, we conclude that $z=1$. Thus, \textbf{Z}$(G)\cap N=1$. $\Box$\\

 Recall that a conjugacy class $A$ of a group $G$ is said to be real if $A^{-1}=A$, and an element $g \in G$ is real if there exists $x \in G$ such that $g^{x}=g^{-1}$. Of course, all elements of
 order 2 in $G$  are real. We note that if $A$ is the unique class of size $|A|$ inside a certain set of classes containing $A^{-1}$, then the fact that $|A|=|A^{-1}|$ proves that $A$ is necessarily
 real. Another elementary property, which is used in the sequel, is that if $x^G$ is real and $|x^{G}|$ is odd, then $x^{2}=1$. \\

{\bf Theorem 4.1.2.} {\it Let $N$ be a normal subgroup of a group $G$ such that $\Gamma_G(N)$ has two vertices and no edge. Then $N$ is a $2$-group or a Frobenius group with  $p$-elementary abelian
kernel $K$, and complement $H$, which is cyclic of order $q$, for two different primes $p$ and $q$. In particular, $|N|=p^n q$ with $n\geq 1$.}\\

{\it Proof.} By Theorem 2.2 we know that either $N=P\times A$ with $A\leq {\rm \textbf{Z}}(G)$ or $N$ is quasi-Frobenius with abelian kernel and complement. In the first case, by Lemma 4.1.1 we obtain
that $A$ is trivial. On the other hand, notice that one of the $G$-classes of the graph has odd size, so it is real, and this implies that $N$ has elements of order 2. Consequently, $N$ is a
$2$-group.

Suppose now that $N$ is quasi-Frobenius with kernel $K$ and complement $H$, both abelian. It is trivial that there must exist elements of order $p$ and order $q$ in $N$ for two different primes $p$
and $q$. If \textbf{Z}$(N)\neq 1$, then $N$ would have elements of at least three different orders. Thus, \textbf{Z}$(N)=1$, which means that $N$ is a Frobenius group. Now, since $p$ and $q$ are the
only possible orders for the elements of $N$, we have that all nontrivial elements in $K$ must have the same order, say $p$, and since $K$ is abelian, we conclude that $K$ is $p$-elementary abelian.
Furthermore, it is well-known that the Sylow subgroups of $H$ are cyclic or generalized quaternion. Nevertheless, since all elements of $H$ have the same order, we conclude that $H$ is cyclic of order
$q$. $\Box$\\

{\bf Example 2.} Let $G=S_{4}$ and let $N=A_{4}$. The $G$-classes of $N$ have sizes $\lbrace 1, 3, 8 \rbrace$ and accordingly, $\Gamma_{G}(N)$ has two vertices and no edge.

\subsection{$\Gamma_{G}(N)$ with two vertices and one edge}
We remark that this case does not occur when dealing with the ordinary graph $\Gamma(G)$. In fact,  as we pointed out in the Introduction, the complete list of those groups $G$ whose graph $\Gamma(G)$
has no triangles proves that the graph $\Gamma(G)$  can never have exactly two vertices joined by an edge.\\

Hereafter $\pi(G)$ will denote the set of primes dividing the order of $G$. We prove first the following property, which will also be used to prove Theorem A.

\bigskip

{\bf Lemma 4.2.1.} {\it Let $N$ be a normal subgroup of a group $G$. Suppose that $\Gamma_G(N)$ is non-empty, has no triangles and that $|N|$ is not a prime power. Then either ${\bf Z}(G)\cap N=1$ or
$|${\rm{\textbf{Z}}}$(G)\cap N|_2= 2$.}\\

{\it Proof.} Suppose that ${\bf Z}(G)\cap N>1$. Since $N$ is not abelian and has not prime power order, we can trivially take a  $q$-element $z \in {\rm{\textbf{Z}}}(G)\cap N$  and a $p$-element $x\in
N \setminus {\bf Z}(G)\cap N$ for two distinct primes $p$ and $q$. If $q>2$, then $x^G$, $(zx)^G$ and $(z^2x)^G$ would be three distinct $G$-classes of $N$ forming a triangle in $\Gamma_G(N)$, a
contradiction.  This implies that $q=2$. Now, if $|{\bf Z}(G)\cap N|_2 > 2$, we can take two distinct nontrivial $2$-elements $z_1, z_2\in {\rm{\textbf{Z}}}(G)\cap N$, whence the $G$-classes $x^G$,
$(xz_1)^G$ and $(xz_2)^G$ are distinct too and again we have a triangle in $\Gamma_G(N)$. Therefore, $|{\bf Z}(G)\cap N|_2 = 2$.  $\Box$\\

{\bf Theorem 4.2.2.} {\it Let $N$ be a normal subgroup of a group $G$ such that $\Gamma_{G}(N)$  has exactly two vertices and one edge. Then one of the following possibilities holds:

\begin{enumerate}
\item $N$ is a $p$-group for a prime $p$.
\item $N=P \times Q$ with $P/({\bf Z}(G)\cap P)$ an elementary abelian $p$-group  with $p$ an odd prime, and $Q \subseteq  {\bf Z}(G)\cap N $ and $Q \cong {\mathbb Z}_{2}$.
\item $N$ is a Frobenius group with $p$-elementary abelian kernel  $K$ and complement $H \cong {\mathbb Z}_{q}$ for some distinct primes  $p$ and $q$. In particular, $|N|=p^{a}q$ for  some $a\geq
    1$ and the $G$-classes of $N$ have cardinality $1$, $(p^{a}-1)$ and $p^{a}(q-1)$.
\end{enumerate} }

{\it Proof.} We divide the proof into two cases, depending on whether $\overline{N}=N /({\rm\textbf{Z}}(G) \cap N)$ is a $p$-group or not.

a) Suppose first that $\overline{N}$ is a $p$-group for some prime $p$. If ${\rm \textbf{Z}}(G) \cap N$ is also a $p$-group, we trivially have the first assertion. Thus, we can assume  that ${\bf
Z}(G) \cap N$ is nontrivial and  apply Lemma 4.2.1, so $|{\bf Z}(G)\cap N|_2=2$.

  Now, observe that $N$ is nilpotent because $\overline{N}$ is a $p$-group. We prove that ${\rm \textbf{Z}}(G) \cap N$ is a $\lbrace 2,p\rbrace$-group and consequently, $N$ will be a $\{2,p\}$-group
  too. In fact, if there exists  $r \in \pi(\textbf{Z}(G) \cap N)$ with $r \neq 2, p$, then we can find in $N \setminus({\rm\textbf{Z}}(G)\cap N)$ elements of order divisible by $2r$ and $pr$, which
  is a contradiction. Therefore, $N$ is a nilpotent $\{2, p\}$-group and we can write $N=P \times Q$ with  $P$ a $p$-group ($p$ odd) and $Q$ a 2-group. Also, we have proved that $Q\cong {\mathbb Z}_2$
  and $Q\subseteq {\bf Z}(G)$. Now, the $G$-classes of $N$ are exactly the products of the $G$-classes of $P$ by the two  (central) classes of $Q$ and, by the structure of the graph, it follows that
  $P$ has exactly one nontrivial $G$-class. This forces $P/({\bf Z}(G)\cap P)$ to be $p$-elementary abelian by Theorem 3.1, and this completes the structure of $N$ described in 2.

b) Assume that  $\overline{N}$ has not prime power order and take two distinct primes $p, q \in \pi(\overline{N})$.  Hence there exist $p$-elements and $q$-elements in $ N \setminus({\rm\textbf{Z}}(G)
\cap N)$. Moreover, all $p$-elements of $ N \setminus ({\rm\textbf{Z}}(G) \cap N)$ form a single non-central $G$-class and the $q$-elements form the other one. No more primes can be involved in
$\pi(\overline{N})$. We note that ${\bf Z}(G)\cap N=1$ because otherwise, the product of one element of prime order of $N\cap {\bf Z}(G)$ by the elements of $ N \setminus ({\rm\textbf{Z}}(G) \cap N)$
produce elements in $N \setminus({\rm\textbf{Z}}(G)\cap N)$ whose order is divisible by more than one prime. This means that $N$ is a $\lbrace p, q\rbrace$-group, and in particular, by
$p^{a}q^{b}$-Burnside's Theorem, $N$ is solvable. Let $M\leq N $ be a minimal normal subgroup of $G$, which is elementary abelian, say for instance, a $p$-group. As $M$ is proper in $N$ and is the
union of $G$-classes, $M$ is equal to the union of the trivial class and one of the non-central $G$-classes of $N$. Then every element of $N \setminus M$ is a $q$-element, so $M \in $ Syl$_{p}(N)$.

Now, if $Q \in $ Syl$_{q}(N)$, then $Q$ acts Frobeniusly on $M$, because $N$ cannot have elements of order divisible by $pq$.
 Moreover, $Q$ must be cyclic or generalized quaternion. As we have seen that all nontrivial $q$-elements have the same order, we conclude that $Q$ is cyclic of order $q$. $\Box$\\

The following examples illustrate the results obtained in the above theorem.

\bigskip
{\bf Example 3.} Let $G=$ SL$(2,3)\times \langle z \rangle$ with $\langle z \rangle \cong {\mathbb Z}_{2}$. The group SL$(2,3)$ has a normal subgroup  $K\cong Q_8$ and we set $N=K\times \langle z
\rangle$. Then ${\rm \textbf{Z}}(G)\cap N = {\rm \textbf{Z}}(K)\times \langle z\rangle$, so $|{\rm \textbf{Z}}(G)\cap N|=4$, and the non-central $G$-classes of $N$ are:
\begin{center}
$(x,1)^{G}=\lbrace (x,1) | x \in K\setminus {\rm\textbf{Z}}(K)\rbrace$ \\
$(x,z)^{G}=\lbrace (x,z) | x \in K\setminus {\rm\textbf{Z}}(K)\rbrace$.
\end{center}
Moreover, $|(x,1)^{G}|=|(x,z)^{G}|=6$, so $\Gamma_{G}(N)$ has two vertices joined by an edge. This example corresponds to case 1.\\

We show two examples of case 2. In the first one, $N$ is abelian whereas in the second, it is not. We also give an example of case 3 of Theorem 4.2.2.\\

{\bf Example 4.} Let $p$ be an odd prime and let
$N=\langle x_{1}\rangle \times \ldots \times \langle x_{n}\rangle \times \langle z \rangle $ with $\langle x_{i} \rangle \cong {\mathbb Z}_{p}$ for $i=1,\ldots, n$ and $\langle z \rangle \cong
{\mathbb Z}_{2}$, so $|N|=2p^{n}$. Let us consider $H\cong $ GL$(n,p)$ the group of automorphisms of $K=\langle x_{1}\rangle \times \ldots \times \langle x_{n}\rangle$ and suppose that $H$ acts
trivially on $\langle z \rangle$. We construct the corresponding semidirect product $G=NH$. Then ${\rm \textbf{Z}}(G)\cap N=\langle z\rangle$, so $|{\rm \textbf{Z}}(G)\cap N|=2$, and the non-central
$G$-classes of $N$ are
\begin{eqnarray*}
\lbrace x_{1}^{i_{1}}\ldots x_{n}^{i_{n}} \mid (i_{1}, \ldots, i_{n})\neq (0, \ldots, 0), 0 \leq i_{j} \leq p-1\rbrace.\\
\lbrace x_{1}^{i_{1}}\ldots x_{n}^{i_{n}}z \mid (i_{1}, \ldots, i_{n})\neq (0, \ldots, 0), 0 \leq i_{j} \leq p-1\rbrace.
\end{eqnarray*}

Thus, $\Gamma_{G}(N)$ has  exactly two vertices, corresponding to the above two classes, both of size $p^{n}-1$ and trivially joined by an edge.  \\

{\bf Example 5.} Let $P$ be the extraspecial $p$-group of order $p^{3}$ and exponent $p$ ($p$ odd):
$$P=\langle a, b, c \mid a^{p}=1,\, b^{p}=1,\, c^{p}=1,\, [a,b]=[b,c]=1, [a,c]=b \rangle.$$
We have $P= \lbrace a^{i}b^{j}c^{k} \mid 0\leq i, j, k \leq p-1\rbrace$ and observe that  ${\bf Z}(P)=\langle b\rangle$.
The group Aut$(P)$ is known (see for instance, Theorem 20.8 of \cite{Doerk}), has order $p^{3}(p-1)^2(p+1)$ and  it turns out that every $\varphi\in$ Aut$(P)$ is determined by: $\varphi(a)=
a^{i}b^{j}c^{k}$; $\varphi(b)= b^{m}$, with $0< m \leq p-1$; and $\varphi(c)= a^{q}b^{r}c^{s}$,
satisfying further that $m\equiv is-qk \not \equiv$ $0\,$ ({\rm mod} $p$) with $0\leq i,\, j,\, k,\, q,\, r,\, s\, \leq p-1$.

Notice that there always exists an automorphism of $P$ leading $a$ to any non-central element of $P$.  Furthermore,  we can set $m=1$, that is, the image of $a$ can be chosen in such a way that all
elements of ${\rm\textbf{Z}}(P)$ are fixed element by element. Set $H=\{\varphi \in {\rm Aut}(P)/ \varphi(b)=b\}$
and let us consider the semidirect product $G=(P \times \langle z\rangle)\rtimes H$, where $\langle z \rangle \cong {\mathbb Z}_{2}$, and $H$ acts naturally on $P$ and trivially on $\langle z\rangle$.
Let $N=P\times \langle z\rangle$. It follows that ${\rm \textbf{Z}}(G)\cap N=\langle b\rangle \times \langle z\rangle$, so $|{\rm \textbf{Z}}(G)\cap N|=2p$, and the non-central $G$-classes of $N$ are
$(x,1)^{G}=\lbrace (x,1)  / x \in P\setminus ${\rm \textbf{Z}}$(P)\rbrace$ and $(x,z)^{G}=\lbrace (x,z) / x \in P\setminus ${\rm \textbf{Z}}$(P)\rbrace$.
Then $|(x,1)^{G}|=|(x,z)^{G}|=p^{3}-p$ and $\Gamma_{G}(N)$ has two vertices joined by an edge.\\

{\bf Example 6.} We use the semilinear affine group $\Gamma(p^n)$ for appropriate $p$ and $n$. We remind the construction of this group, although for a more detailed description we refer the reader to
Section 2 of \cite{Manz}. Let $GF(p^{n})$ be the finite field of $p^{n}$ elements.  The multiplicative group $H=GF(p^{n})^{*}$ is cyclic of order $p^{n}-1$ and acts on the additive group of
$GF(p^{n})$, say $K=\langle x_{1}\rangle \oplus \ldots \oplus \langle x_{n}\rangle \cong {\mathbb Z}_{p}\times \ldots \times {\mathbb Z}_{p}$ in the following way: if $a \in H$ and $x \in K$, then
$x^{a}=ax$. This action is Frobenius, so $KH$ is a Frobenius group with abelian kernel and complement. On the other hand, we know that $\alpha$, defined by $x^{\alpha}=x^{p}$ for all $x \in K$, is
another automorphism of $K$ of order $n$ in such a way that $H\langle \alpha \rangle\leq $ Aut$(K)$. Then $\Gamma(p^n)$ is defined by the semidirect product $\Gamma(p^n)=K\rtimes (H\langle \alpha
\rangle)$.

Now, suppose that $p$ is odd, $n=2$ and $S\leq H$ is  a cyclic subgroup of order $s=3$ (note that we need to assume that 3 divides $p^{2}-1$). One can easily check that $N=KS$ is normal in $G$, and
that $N$ is a Frobenius group with abelian kernel and complement. The two nontrivial $G$-classes of $N$ consist of the $p^{2}-1$ elements of $K\setminus \lbrace 1 \rbrace$, and all the elements of
$\bigcup_{g \in G}(S^{g}\setminus \lbrace 1 \rbrace )=\bigcup_{k \in K}(S^{k}\setminus
\lbrace 1 \rbrace)$, which are exactly $(|S|-1)|K|=2p^{2}$ elements. As $(2p^{2}, p^{2}-1)=2$,  we obtain case 3.

\section{$\Gamma_{G}(N)$ with three vertices}

We distinguish the following three possibilities: $\Gamma_{G}(N)$ is disconnected with three vertices, $\Gamma_{G}(N)$ is connected having exactly one triangle or $\Gamma_{G}(N)$ has three
vertices in a line. The former two  possibilities may happen for the ordinary graph $\Gamma(G)$. In fact, $\Gamma(G)$ is disconnected with three vertices if and only if $G\cong D_{10}$ or $G\cong A_4$
(see \cite{FangZhang}). In subsection 5.2, we will study the case in which $\Gamma(G)$ is exactly a triangle, by showing that $G\cong Q_8$ or $G\cong D_8$. However, the latter case, three vertices in a line, do not occur for the ordinary graph (\cite{FangZhang}).

\subsection{$\Gamma_{G}(N)$ disconnected with three vertices}

{\bf Theorem 5.1.1.} {\it Let $N$ be a normal subgroup of a group $G$. If $\Gamma_{G}(N)$ has three vertices and one edge, then  $N$ is a $\{p,q\}$-group for two primes $p$ and $q$. Furthermore,
either
\begin{enumerate}

\item[1)]  $N$ is a $p$-group, or

\item[2)] $N$ is a quasi-Frobenius group with abelian kernel and complement. In this case, $|N \cap$ {\rm\textbf{Z}}$(G)| =1$ or $2$.
\end{enumerate}
}

{\it Proof.}  The graph $\Gamma_{G}(N)$ has the form:
\begin{center}
\begin{tikzpicture}
  [scale=.5,auto=left,every node/.style={circle,fill=black!20, scale=0.7}]
  \node (n1) at (0,0) {};
  \node (n2) at (3,0)  {};
  \node (n3) at (6,0)  {};
\foreach \from/\to in {n1/n2,n3}
\draw (\from) -- (\to);
\end{tikzpicture}
\end{center}
 and by Theorem 2.2, we know that either $N$ is quasi-Frobenius with abelian kernel and complement, or $N=P\times A$ where $P$ is a $p$-group and $A\leqslant {\rm \textbf{Z}}(G)$.  Suppose the second
 case. We can certainly assume that $A$ is a $p'$-group. If $A>1$, then every non-central $G$-class of $N$ is the product of a non-central $G$-class of $P$ by  a (central) class of $A$. By hypothesis,
 the number of non-central $G$-classes of $N$ is 3 and this occurs if and only if either $|A|=3$ and $P$ has exactly one non-central $G$-class, or $|A|= 1$ and $P$ has 3 non-central $G$-classes.
 However, in the former case, the three $G$-classes would form a triangle, so we conclude that $N$ is a $p$-group.

 Assume now that $N$ is quasi-Frobenius with abelian kernel and complement. We claim that $N$ is a $\{p,q\}$-group for two primes $p$ and $q$. First, suppose that $\textbf{\mbox{Z}}(G) \cap N=1$. If
 $|N|$ is divisible by three distinct primes, then the existence of exactly three $G$-classes implies that all elements of $N$  necessarily have prime order and, by Theorem 2.5, $N$ is a
 $\{p,q\}$-group or $N\cong A_5$. Since $N$ is clearly solvable, we are finished.  Now, if ${\bf Z}(G)\cap N> 1$, by Lemma 4.2.1, we can assume $|\textbf{\mbox{Z}}(G)\cap N|_2=2$. If
 $|N/(\textbf{\mbox{Z}}(G) \cap N)|$ is divisible by two primes $p$ and $q$, both different from 2, then we can multiply a central 2-element by  a $p$-element and a $q$-element of $N$ in order to get
 non-central elements of order divisible by $2p$ and $2q,$ respectively, and this contradicts the hypothesis. Hence, for instance $p=2$, and  $N$ is a $\{2,q\}$-group. In order to conclude this case
 it is enough to see that $|N \cap {\rm\textbf{Z}}(G)|_q=1$. In fact, if there exists a nontrivial $q$-element in $x\in {\bf Z}(G)\cap N$, we choose a $2$-element $y\in N\setminus ({\bf Z}(G)\cap N)$,
 and then $y^G$, $(xy)^G$ and $(x^2y)^G$ are three non-central $G$-classes of $N$ with the same cardinality, a contradiction.  $\Box$\\

{\bf Example 7.} If we consider $G=N=A_{4}$, then $\Gamma_{G}(N)$ has three vertices and one edge. This example illustrates case 1.
On the other hand, the group of the Small groups library of GAP with number Id (324,
8) has an abelian normal subgroup $N \cong {\mathbb Z}_3 \times {\mathbb Z}_3$ whose $G$-class sizes are 1, 2, 3 and 3. This example shows case 2.

\subsection{$\Gamma_{G}(N)$ with exactly one triangle}

In this section we prove Theorems B and C. We analyze first the ordinary case.

\bigskip

{\it Proof of Theorem B.} By the class equation, we have $$|G|=|{\bf Z}(G)| + |a^G|+ |b^G| +|c^G|,$$ where $a^G, b^G$ and $c^G$ are the three non-central classes of $G$.  Since ${\bf Z}(G) \subset
{\bf C}_G(x) \subset G$, for $x = a,b,c$, then
$$1=\frac{|{\bf Z}(G)|}{|G|} + \frac{1}{|{\bf C}_G(a)|}+ \frac{1}{|{\bf C}_G(b)|} +\frac{1}{|{\bf C}_G(c)|} \leq \frac{|{\bf Z}(G)|}{|G|}+\frac{3}{2|{\bf Z}(G)|}.$$
Therefore, $$\frac{3}{4} \leq 1 - \frac{|{\bf Z}(G)|}{|G|} \leq \frac{3}{2|{\bf Z}(G)|}$$ and thus, $|{\bf Z}(G)| \leq 2$. If $|{\bf Z}(G)| =2$, then $$1 \leq \frac{|{\bf Z}(G)|}{|G|} +
\frac{3}{2|{\bf Z}(G)|} \leq \frac{2}{|G|} + \frac{3}{4}.$$ This implies that $|G| \leq 8$ and necessarily $G\cong Q_8$ or $G\cong D_8$.

Suppose now that ${\bf Z}(G)=1$ and we seek a contradiction. If $|G|$ is divisible by three different primes, then by applying the hypotheses, $G$ should have all nontrivial elements of prime
order, so by applying Theorem 2.5, we get $G\cong A_5$, which has exactly four non-central classes, a contradiction. Thus, we can assume that $G$ is a $p$-group or a $\{p,q\}$-group for two distinct
primes $p$ and $q$. Since $\Gamma(G)$ consists exactly of one triangle, one prime, say for instance $p$, must divide the size of every non-central class of $G$, and then by the class equation, $p$
must divide $|{\bf Z}(G)|$, a contradiction too, as wanted.  $\Box$

\bigskip

{\it Proof of Theorem C.} If $N$ is a $p$-group, we have case 1, so we can assume that $|N|$ is divisible by at least two distinct primes. Let $\overline{N}=N/ ({\rm\textbf{Z}}(G)\cap N)$ and
distinguish two possibilities, whether $\overline{N}$ is a CP-group or not.\\

a) $\overline{N}$ is not a CP-group. Let $\overline{x}\in \overline{N}$ of order $pq$ with $p\neq q$. We write $x=x_px_q$ and may assume that $x_p$ and $x_q$ are non-central  of order $p$ and $q$,
respectively, so the three non-central $G$-conjugacy classes of $N$ are $x_p^G$, $x_q^G$ and $x^G$. Notice that ${\rm\textbf{Z}}(G)\cap N=1$, otherwise we could obtain more than three non-central
$G$-classes in $N$. As a result, $N$ is a $\{p,q\}$-group. We prove now that ${\bf Z}(N)>1$. Suppose on the contrary that ${\rm \textbf{Z}}(N)=1$ and let $P$ and $Q$ be a Sylow $p$ and $q$-subgroup of
$N$, respectively. Let $z \in {\rm \textbf{Z}}(P)$. Since all $p$-elements of $N$ are $G$-conjugate, they lie in $x_{p}^{G}$ and we have $|z^{G}|=|x_{p}^{G}|=|G:{\rm
\textbf{C}}_{G}(x_{p})N||x_{p}^{N}|$. Moreover, ${\rm \textbf{C}}_{G}(z)N=({\rm \textbf{C}}_{G}(x_{p})N)^{g}$, so  $|z^{N}|=|x_{p}^{N}|$ is a power of $q$. Consequently, $|z^{N}|=|x_{p}^{N}|=|n^{N}|$
for every $p$-element $n \in N$. Similarly, we get that $|x_{q}^{N}|$ is power of $p$ and $|x_{q}^{N}|=|m^{N}|$ for every $q$-element $m\in N$. Furthermore, as the $N$-class sizes of $x_p$ and $x_q$
are coprime, we have $N={\rm \textbf{C}}_{N}(x_{p}){\rm \textbf{C}}_{N}(x_{q})$, and since ${\bf C}_N(x)={\bf C}_N(x_p)\cap {\bf C}_N(x_q)$, it follows that $|x^{N}|=|x_{p}^{N}||x_{q}^{N}|$. Thus, the
set of class sizes of $N$ is cs$(N)=\lbrace 1, |x_{p}^{N}|, |x_{q}^{N}|, |x^{N}|\rbrace$ and we can apply the main result of \cite{Camina} to deduce that $N$ is nilpotent. In particular, it would have
nontrivial center, a contradiction. This contradictions proves that ${\rm \textbf{Z}}(N)\neq 1$, as wanted. Now, if  ${\rm \textbf{Z}}(N)$ has $p$-elements (analogously for $q$), then $P\subseteq {\rm
\textbf{Z}}(N)$ and $N=P\times Q$. Since the three non-central classes of $N$ are obtained as the products of classes of the factors, we deduce  that both $P$ and $Q$ can only have a non-central
$G$-conjugacy class. By Theorem 3.1, $P$ and $Q$ are elementary abelian. Therefore, we get case 2.\\

b) $\overline{N}$ is a CP-group. By assumption, we have that $|\overline{N}|$ may be divisible at most by three different primes. We distinguish three subcases:\\

 b.1)  $\overline{N}$ is a $p$-group. Let us see that $\overline{N}$ has exponent $p$. Let $\overline{x}\in \overline{N}$ with $x\in N$ a $p$-element and suppose that $o(x)=p^n$ with $n>1$. We know
 that there exists a $q$-element $z \in {\rm\textbf{Z}}(G)\cap N$  with $q\neq p$ because $N$ is not a $p$-group. In fact, $N$ has a central Sylow $q$-subgroup and we can write $N=P\times Q$. Thus,
 $x^{G}$, $(x^{p})^{G}$, $(zx)^{G}$ and $(zx^{p})^{G}$ would be four different conjugacy classes, a contradiction. As a consequence, $o(x)=p$ and $o(\overline{x})=p$. Now, let us see that
 $|{\rm\textbf{Z}}(G)\cap N|$ is divisible at most by $p$ and $q$ with $q\neq p$. If another prime $r$ divides $|{\rm\textbf{Z}}(G)\cap N|$, then $N$ has non-central $p$-elements, $\lbrace p,
 q\rbrace$-elements, $\lbrace p, r\rbrace$-elements and $\lbrace p, q, r\rbrace$-elements, again a contradiction. We prove now that $Q\cong {\mathbb Z}_{3}$. Since $N=P\times Q$, we have that the
 $G$-classes of $N$ are the product of the $G$-classes of $P$ and $Q$. Particularly, since $Q\leq {\bf Z}(G)$,  the non-central $G$-classes of $N$ are the product of the $G$-classes of $P$ (distinct
 from 1) multiplied by all $G$-classes of $Q>1$. As there are three non-central $G$-classes in $N$, then $Q$ necessarily has exactly three central $G$-classes, so $Q\cong {\mathbb Z}_{3}$ and we get
 case 3.\\

b.2)  $\overline{N}$ is $\lbrace p, q\rbrace$-group. Note that $N$ is also a $\lbrace p, q\rbrace$-group. On the contrary, if there is a prime $r\neq q,p$ dividing $|{\rm\textbf{Z}}(G)\cap N|$, then
there are non-central elements whose order is divisible by $pr$ and $qr$ in $N$, giving a contradiction. We distinguish two cases. If ${\rm\textbf{Z}}(G)\cap N=1$, then $N$ is a CP-group. By applying
Theorem 2.4 and taking into account that there are only three non-central $G$-classes in $N$, one easily obtains the possibilities described in case 5. Suppose now that ${\rm\textbf{Z}}(G)\cap N\neq
1$ and we see that $|{\rm\textbf{Z}}(G)\cap N|=2$. Suppose that for one prime, say $q$ for instance, $|{\bf Z}(G)\cap N|_q>2$ and take two distinct nontrivial $q$-elements  $z_{1}, z_{2}\in
{\rm\textbf{Z}}(G)\cap N$.  We choose  a $p$-element $x\in N\setminus {\rm\textbf{Z}}(G)\cap N$ and we have that $x^G$, $(z_{1}x)^{G}$ and $(z_{2}x)^{G}$ are three distinct $G$-classes. However, we
also have non-central $q$-elements in $N$, so there are at least four non-central $G$-classes inside $N$, a contradiction. Hence $|{\rm\textbf{Z}}(G)\cap N|=2$, and we obtain case 4. \\

b.3) $\overline{N}$ is $\lbrace p, q, r\rbrace$-group. It clearly follows that ${\rm\textbf{Z}}(G)\cap N=1$ by counting order of elements. This means that $N$ is $\lbrace p, q, r\rbrace$-group and
that all elements have prime order. By Theorem 2.5, we have $N\cong A_{5}$. Now, it is known that $G/{\bf C}_G(N)$ is immersed in Aut$(N)\cong S_5$. Then

$$N\cong N{\bf C}_G(N)/{\bf C}_G(N)\subseteq G/{\bf C}_G(N)\cong A_5 \quad  {\rm  or } \quad S_5.$$
We have two possibilities: either $G =N{\bf C}_G(N)$ or $G/{\bf C}_G(N)\cong S_{5}$. In the first case, $G=N\times {\bf C}_G(N)$ and the $N$-classes of $N$ are $G$-classes, so $N$ has four non-central
$G$-classes, a contradiction. Consequently, $G/{\bf C}_G(N) \cong S_{5}$. It is enough to take $K={\bf C}_G(N)$ and $x\in G\setminus (N\times K)$, so as to obtain case 6. $\Box$\\

{ \bf Example 8.} We show examples for each of the cases of Theorem C. The groups $G=D_8$ and $G=Q_8$ with $N=G$ are examples of case 1.

Let us take $A\cong {\mathbb Z}_6$ acting Frobeniusly on  a group $H\cong {\mathbb Z}_7$. Let us consider the action of $A$ on another group $K\cong {\mathbb Z}_3$ in such a way that the kernel of the
action  is $A_0\cong {\mathbb Z}_3$ and $A/A_0\cong {\mathbb Z}_2$ acts Frobeniusly on $K$. Then we construct the corresponding semidirect product $G=(H\times K)\rtimes A$ and  take $N=K\times H$. We
have that the sizes of the $G$-classes  of $N$ are 1, 2, 6 and 12. This example shows case 2.

Let $G=$ SL$(2,3)\times {\mathbb Z}_3$ and let $N=K\times {\mathbb Z}_3$, where $K$ is the normal subgroup of SL(2,3) isomorphic to $Q_8$. Then the three non-central $G$-classes of $N$ have size  6,
so  we have an example of case 3.

 If we consider $G=$ GL$(2,3)$ with $N=$ SL$(2,3)$, then the $G$-classes of $N$ have size 1, 1, 6, 8, 8, and this illustrates case 4.

 The same construction of Example 6, that is,  $G=\Gamma(p^2)$ may work for case 5. We have to take  $s= 5$ instead of $3$, (so we need that  5 divides $p^2 -1$),  and set $N=KS$ with $S\leq H$ a
 cyclic group of order $5$. Then the sizes of the $G$-classes of $N$ are $\{1, p^2-1, 2p^2, 2p^2\}$, so $\Gamma_G(N)$ is just a triangle when $p$ is odd.

Finally, of course $G=S_5$ with $N=A_5$ is an example of case 6. \\

{\bf Remark.} As we have pointed out in the introduction,
the $S_3$-conjecture does not hold for conjugacy classes lying in a normal subgroup. The examples given in Sections 3, 4 and 5 show that even there are infinitely many normal subgroups having all
$G$-classes (central and non-central) with distinct size.

\subsection{$\Gamma_{G}(N)$ with three vertices connected in a  line}
{\bf Theorem 5.3.1.} {\it Let $N$ be a normal subgroup of a group $G$. If $\Gamma_{G}(N)$ has three vertices in a line then {\rm \textbf{Z}}$(G)\cap N=1$ and one of  the following cases is satisfied:
 \begin{enumerate}
\item $N$ is  a $2$-group of exponent at most 4.
\item $N= P \times Q$, where $P$ and $Q$ are elementary abelian $p$ and $q$-groups.
\item $N$ is a Frobenius group with complement isomorphic to ${\mathbb Z}_q$, ${\mathbb Z}_{q^2}$ or $Q_8$. In the former case, the kernel of $N$ is a $p$-group with exponent $\leq p^2$ and in the
    two latter cases, the kernel of $N$ is $p$-elementary abelian.
\end{enumerate}}
In all cases, $|N|$ is divisible by at most two prime numbers.\\

{\it Proof.} The graph $\Gamma_{G}(N)$ has the form:
 \begin{center}
\begin{tikzpicture}
  [scale=.5,auto=left,every node/.style={circle,fill=black!20, scale=0.7}]
  \node (n1) at (0,0) {};
  \node (n2) at (3,0)  {};
  \node (n3) at (6,0)  {};
\foreach \from/\to in {n1/n2,n2/n3}
\draw (\from) -- (\to);
\end{tikzpicture}
\end{center}

Since the vertices of the ends, say $x^{G}$ and $y^{G}$, are not joined, then  $x^G$ and $y^G$ are the only $G$-classes with cardinality $|x^G|$ and $|y^G|$, respectively. If $z \in {\rm
\textbf{Z}}(G)\cap N$, then $(zx)^G = x^G$ and $(zy)^G = y^G$, so $|\langle z \rangle|$ divides $|x^G|$ and $|y^G|$. Therefore, ${\rm \textbf{Z}}(G)\cap N=1$. We will distinguish three cases:\\

a) Suppose that $N$ is a $p$-group. By the graph form $x^G$ and $y^G$ are real classes and obviously one of them needs to have odd size, say $y^G$. Thus, the order of $y$ is 2 and $p=2$. Notice that
necessarily the order of $x$ is also 2, otherwise  $|(x^2)^G|$ would divide $|x^G|$, and this would force $\Gamma_G(N)$ to have a triangle. If $z^G$ is the third $G$-class of $N$, then the order of
$z$ is at most 4 and we obtain case 1.\\

b) Suppose now that $N$ is not a  $p$-group and not a CP-group. Hence it has elements whose order is divisible by two different primes $p$ and $q$. Let $x=x_{p}x_{q}$ be one of them, where $x_{p}$ and
$x_{q}$ are the (non-central) $p$-part and $q$-part, respectively. Hence, the three $G$-classes of $N$ are $x_{p}^{G}$, $x_{q}^{G}$ and $x^{G}$, so $N$ is a $\{p,q\}$-group. Let us prove that ${\rm
\textbf{Z}}(N)\neq 1$. Suppose that ${\rm \textbf{Z}}(N) = 1$. Since all $p$-elements of $N$ lie in $x_p^G$, any $p$-element $w\in N$ satisfies $w = x_p^g$, for some $g \in G$ and
$$|x_{p}^{G}|=|G:{\bf C}_{G}(x_{p})N| |x_{p}^{N}|= |G:{\bf C}_{G}(w)N| |w^{N}|= |w^G|. $$
  It follows that $|x^N| = |w^N|$. We can argue similarly for the $N$-classes of $q$-elements and $\lbrace p, q \rbrace$-elements, which have sizes $|x_{q}^{N}|$ and $|x^{N}|$, respectively. In
  addition, since $|x_{p}^{N}|$ and $|x_{q}^{N}|$ divide $|x^{G}|$ and $|x_q^G|$ respectively, we have $(|x_{p}^{N}|, |x_{q}^{N}|)=1$ by the structure of the graph. Furthermore, $|x_{p}^{N}|$ and
  $|x_{q}^{N}|$ must be a power of $p$ and a power of $q$ greater than $1$ because ${\rm \textbf{Z}}(N)=1$. On the other hand, $N={\rm \textbf{C}}_{N}(x_{p}){\rm \textbf{C}}_{N}(x_{q})$, ${\rm
  \textbf{C}}_{N}(x)={\rm \textbf{C}}_{N}(x_{p})\cap {\rm \textbf{C}}_{N}(x_{q})$ and this leads to $|x^{N}|=|x_{p}^{N}||x_{q}^{N}|$. Therefore, the set of class sizes of $N$ is ${\rm cs}(N)=\lbrace
  1, |x_{p}^{N}|, |x_{q}^{N}|, |x^{N}|\rbrace$ and by the main result of \cite{Camina}, $N$ is nilpotent and has nontrivial center, a contradiction. Now, let us prove that $N$ is abelian. If $p$
  divides $|{\rm \textbf{Z}}(N)|$, since all $p$-elements of $N$ are $G$-conjugate, $N$ has a central Sylow $p$-subgroup, so $N$ is nilpotent and factorize $N=P\times Q$. As there are only $3$
  nontrivial $G$-classes, this implies that $P$ and $Q$ are minimal normal subgroups of $G$, so $P$ and $Q$ are elementary abelian. We get case 2.\\

c) Assume that $N$ is a CP-group but not a $p$-group. By Theorem 2.4, we know that $|N|$ is divisible by at most two primes.  We will prove that case 3 of Theorem 4 cannot happen and, consequently, we
obtain the properties described in case 3 of the statement. Suppose that $N/{\rm \textbf{O}}_{p}(N)$ is a $\{p, q\}$-group having cyclic Sylow subgroups and let $K/{\rm \textbf{O}}_{p}(N)$ and $H/{\rm
\textbf{O}}_{p}(N)$  be the kernel and complement of the Frobenius group $N/{\rm \textbf{O}}_{p}(N)$, respectively. Observe that $K/{\rm \textbf{O}}_{p}(N)$ is a cyclic $q$-group and $H/{\rm
\textbf{O}}_{p}(N)$ is a cyclic $p$-group. By the hypotheses, we can write
$${\rm \textbf{O}}_{p}(N)=1\cup x^{G},$$ $$K=1\cup x^{G}\cup y^{G},$$ $$N=1\cup x^{G}\cup y^{G}\cup z^{G},$$
 where $x$ has order $p$ and $y$ has order $q$. Notice that $y^{G}$ is the only class of elements of order exactly $q$ because $K/{\rm \textbf{O}}_{p}(N)\cong {\mathbb Z}_{q}$. Also, we can write
 $K={\rm \textbf{O}}_{p}(N)Q$ where $Q\cong {\mathbb Z}_{q}$ is a Sylow $q$-subgroup of $G$, and $|K|=p^{a}q$. On the other hand, the elements of $N\setminus K$ form a $G$-class, say $z^{G}$.
 Moreover,
 $$\frac{N}{K}\cong \frac{N/{\rm \textbf{O}}_{p}(N)}{K/{\rm \textbf{O}}_{p}(N)}\cong H/{\rm \textbf{O}}_{p}(N)$$
 and $|N|=p^{a+1}q$.
 Therefore,
 $$|x^{G}|= |{\rm \textbf{O}}_{p}(N)|-1=p^{a}-1,$$
  $$|y^{G}|=|K|-|{\rm \textbf{O}}_{p}(N)|=p^{a}(q-1),$$
   $$|z^{G}|=|N|-1-|x^G|-|y^G|= p^{a}q(p-1).$$

    If both $p$ and $q$ are odd, we have that $|x^{G}|$, $|y^{G}|$ and $|z^{G}|$ are even, a contradiction. If $q=2$, then  $|K/{\rm \textbf{O}}_{p}(N)|=2$, and $K/{\rm \textbf{O}}_{p}(N)$ would be
    central in $N/{\bf O}_p(N)$, again a contradiction. Thus, $p=2$ and $|N|=2^{a+1}q$. In particular, $|H/{\rm \textbf{O}}_{2}(N)|=2$. By the Frattini argument, $N=K{\rm \textbf{N}}_{N}(Q)={\rm
    \textbf{O}}_{2}(N){\rm \textbf{N}}_{N}(Q)$. On the other hand,  $[{\rm \textbf{N}}_{N}(Q)\cap {\rm \textbf{O}}_{2}(N), Q ]\subseteq Q \cap {\rm \textbf{O}}_{2}(N)=1$ and, in particular, ${\rm
    \textbf{N}}_{N}(Q)\cap {\rm \textbf{O}}_{2}(N)\subseteq {\rm \textbf{C}}_{N}(Q)\cap {\rm \textbf{O}}_{2}(N)=1$. Then, $|{\rm \textbf{N}}_{N}(Q)|=|N|/|{\rm \textbf{O}}_{2}(N)|=2q$ and there must
    exist a $2$-element $w \in {\rm \textbf{N}}_{N}(Q)$. If $H$ were abelian, then $H\subseteq {\rm \textbf{C}}_{N}({\rm \textbf{O}}_{2}(N))$. However, we know that ${\bf O}_{2'}(N)=1$, so  ${\rm
    \textbf{C}}_{N}({\rm \textbf{O}}_{2}(N))\subseteq {\rm \textbf{O}}_{2}(N)$, which provides a contradiction. As a result, $H$ is not abelian and $z^{G}$ consists of elements of order $4$.
  In particular, $w$ has order $4$ and satisfies $(w{\rm \textbf{O}}_{2}(N))^{2}={\rm \textbf{O}}_{2}(N)$, which means that $w^{2}\in {\rm \textbf{O}}_{2}(N)$.  Then $w^{2}\in {\rm
  \textbf{N}}_{N}(Q)\cap {\rm \textbf{O}}_{2}(N) = 1$, the final contradiction. $\Box$\\

{\bf Example 10.} We take $G=(({\mathbb Z}_2\times {\mathbb Z}_2\times {\mathbb Z}_2)\rtimes {\mathbb Z}_7)\times A_4$ which has been obtained from \cite{gap} by using the code $SmallGroup(672,1258)$.
Then $N={\mathbb Z}_2\times {\mathbb Z}_2 \times {\mathbb Z}_2 \times {\mathbb Z}_2 \times {\mathbb Z}_2$ is a normal subgroup of $G$ and the sizes of the $G$-classes are $\{ 1, 3, 7, 21\rbrace$.
  Thus, $\Gamma_{G}(N)$ has three vertices in a line and $N$ is a $2$-group.	\\

{\bf Example 11.} We consider $N={\mathbb Z}_2\times {\mathbb Z}_2\times {\mathbb Z}_2 \times {\mathbb Z}_3$ and construct an action of $H={\mathbb Z} _7\times {\mathbb Z}_2$ on $K$. We know that
${\mathbb Z}_7$ can act Frobeniusly on ${\mathbb Z}_2\times {\mathbb Z}_2\times {\mathbb Z}_2$, because this group has a fixed-point-free  automorphism of order 7. Also, we consider the trivial action of
${\mathbb Z}_7$ on ${\mathbb Z}_3$. On the other hand, we consider ${\mathbb Z}_2$ acting trivially on ${\mathbb Z}_{2}\times {\mathbb Z}_2\times {\mathbb Z}_{2}$ and Frobeniusly on ${\mathbb Z}_3$.
We set the corresponding semidirect product $G=N \rtimes H$ and we have that $N$ has three nontrivial $G$-classes whose sizes are $2, 7$ and $14$. Accordingly, $\Gamma_{G}(N)$ has three vertices in
a line.\\

{\bf Example 12.} We take $G=(({\mathbb Z}_5\times {\mathbb Z}_5)\rtimes Q_8)\rtimes {\mathbb Z}_3$, which has been obtained with \cite{gap} by using the code $SmallGroup(600,150)$.  Then
$N=({\mathbb Z}_5\times {\mathbb Z}_5)\rtimes Q_8$ is a normal subgroup of $G$ and the sizes of the $G$-classes of $N$ are  $\lbrace 1, 24, 25, 150\rbrace$.  Notice that $\Gamma_{G}(N)$ has three
vertices in a line and $N$ is a Frobenius group with complement $Q_8$.

\section{$\Gamma_{G}(N)$ without triangles}

{\bf Theorem 6.1.} {\it Let $N$ be a normal subgroup of a group $G$ such that $\Gamma_{G}(N)$ has no triangles. Then $N$ is solvable.}\\

{\it Proof.} The case in which $\Gamma_{G}(N)$ is empty is trivial, since $N$ is clearly abelian.  We use induction on $|N|$. Take $N/K$ a chief factor of $N$ with $K<N$ and suppose that $K \neq 1$.
If $\Gamma_G(K)$ has triangles, then  $\Gamma_G(N)$ trivially has triangles too. By induction, we conclude that $K$ is solvable. On the other hand, if $\Gamma_{G/K}(N/K)$ has a triangle, formed by
$(xK)^{G/K}$, $(yK)^{G/K}$ and $(zK)^{G/K}$, then we have that $x^G$, $y^G$ and $z^G$ form a triangle of $\Gamma_G(N)$. In fact, if for instance, $x^G = y^G$, then  $(xK)^{G/K} = (yK)^{G/K}$. Again by
induction, $N/K$ is solvable too. Therefore, we can assume that $N$ is a minimal normal subgroup of $G$. We assume that it is the direct product of isomorphic copies of a non-abelian simple group
$N_1$, otherwise $N$ is solvable and we are finished. Bertram, Herzog y Mann proved in \cite{BerHerMann} that the graph of a simple group is always complete. Thus, we can choose a triangle of
$\Gamma(N_1)$, formed by $a^{N_1}$, $b^{N_1}$ and $c^{N_1}$, such that the orders of $a$, $b$ and $c$ are different prime numbers. Since $N_1\unlhd N \unlhd G$, we deduce that $|n^{N_1}|$ divides
$|n^G|$, for every element $n \in N_1$, so wee obtain  three elements of different orders in $N$, whose corresponding $G$-classes are distinct and form a triangle in $\Gamma_{G}(N)$, a contradiction.
Therefore, $N$ is solvable. $\Box$\\

{\bf Theorem 6.2.} {\it If $\Gamma_{G}(N)$ has no triangles and $N$ is a {\rm CP}-group, then $N$ is a $p$-group for some prime $p$ or a Frobenius group.}\\

{\it Proof. } The case in which $\Gamma_G(N)$ is empty and $N$ is a CP-group, trivially leads to that $N$ is a $p$-group. Suppose that $N$ has not prime power order. By hypothesis, it is clear that
${\rm \textbf{Z}}(N) = 1$ (and ${\rm \textbf{Z}}(G) \cap N = 1$). By Theorem 6.1, we know that $N$ is solvable, and  one of the three options of Theorem 2.4 holds. We prove that case 3 cannot occur.
Suppose that $N$ satisfies case 3 of Theorem 2.4, so $N$ has the following normal series: $$1 <{\rm \textbf{O}}_p(N)   <  K :=  {\rm \textbf{O}}_p(N)Q    < N,$$
where $Q$ is  a Sylow $q$-subgroup of $N$.

Since $K= {\rm \textbf{O}}_p(N)Q$ is a Frobenius group, the $K$-classes of ${\rm \textbf{O}}_p(N) \setminus 1$ have cardinality divisible by $q$, then the cardinalities of the respective $N$-classes
and $G$-classes also are. On the other hand, the elements of $K\setminus {\rm \textbf{O}}_p(N)$ have a $K$-conjugacy class of cardinality divisible by $p$, and again the respective $N$-class sizes and
$G$-class sizes are divisible by $p$ too. Furthermore, if $z$ is a $p$-element of $N \setminus K$, we see that $|z^N|$ is divisible by $p$ and $q$.  It is obvious that $q$ divides $|z^N|$, otherwise
there are elements of order $pq$ in $N$. In addition, $|z^N|$ is divisible by $p$, otherwise ${\bf C}_N(z)$ contains a Sylow p-subgroup of $N$, so in particular, contains ${\rm \textbf{O}}_p(N)$.
Then,  $z\in {\bf C}_N({\rm \textbf{O}}_p(N))\subseteq {\rm \textbf{O}}_p(N)$. As a result, there exist at least three nontrivial $G$-classes in $N$, one of them in ${\rm \textbf{O}}_p(N)$, whose
cardinality is divisible by $q$, another one in $K\setminus {\rm \textbf{O}}_p(N)$ of cardinality divisible by $p$, and a third class in $N\setminus K$ whose cardinality is divisible by $pq$. Since
the $G$-classes are union of $N$-classes, there are not no more nontrivial $G$-classes in $N$, otherwise, there would be a triangle in $\Gamma_G(N)$. Therefore, $\Gamma(N)$ has exactly three vertices
connected in a line. However, this is not possible for the ordinary graph by the main result in \cite{FangZhang}. In conclusion, case 3 of Theorem 2.4 does not happen and so, $N$ is a $p$-group or a
Frobenius group. $\Box$\\

We are ready to prove Theorem A of the Introduction.\\

{\it Proof of Theorem A}. By Theorem 6.1 we know that $N$ is solvable. We will distinguish two cases depending on whether $N/({\rm \textbf{Z}}(G)\cap N)$  is a CP-group or not.

a) Suppose first that $N/({\rm \textbf{Z}}(G)\cap N)$ is not a CP-group. Let $x({\rm \textbf{Z}}(G)\cap N)$ of order divisible by exactly two different prime numbers $p$ and $q$. We can write
$x=x_{p}x_{q}$, where $x_{p}$ and $x_{q}$ are the $p$-part and $q$-part of $x$, respectively, which are trivially non-central in $G$. Since ${\rm \textbf{C}}_{G}(x)={\rm \textbf{C}}_{G}(x_{p})\cap
{\rm \textbf{C}}_{G}(x_q)$, we have that both $|x^{G}_{p}|$ and $|x_{q}^{G}|$ divide $|x^{G}|$, so $x^G$ is connected to $x_{p}^{G}$ and $x_{q}^{G}$. Moreover, by the orders of the elements, these
classes are different. Thus, we have a subgraph of $\Gamma_G(N)$ with the following form:

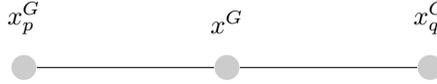
\begin{center}
\begin{tikzpicture}
  [scale=.9,auto=left,every node/.style={circle,fill=black!20}]
  \node [label=$x_{p}^{G}$](n1) at (0,0) {};
  \node [label=$x^{G}$](n2) at (3,0)  {};
  \node [label=$x_{q}^{G}$](n3) at (6,0)  {};
\foreach \from/\to in {n1/n2,n2/n3}
\draw (\from) -- (\to);
\end{tikzpicture}
\end{center}

As $\Gamma_{G}(N)$ has no triangles, then $|x_{p}^{G}|$ and $|x_{q}^{G}|$ are coprime and as a consequence, $G={\rm \textbf{C}}_{G}(x_p){\rm \textbf{C}}_{G}(x_q)$ and $|x^{G}|=|x_{p}^{G}|
|x_{q}^{G}|$. Hence, any prime dividing $|x_{p}^{G}|$ or $|x_{q}^{G}|$ is a divisor of $|x^{G}|$ and consequently, if there exists a vertex in $\Gamma_G(N)$ different from $x^G$, $x_p^G$ and $x_q^G$,
it must lie in another connected component of $\Gamma_{G}(N)$. This means that $\Gamma_{G}(N)$ is disconnected and we get a contradiction by applying Theorem 2.3. Therefore, $\Gamma_{G}(N)$ coincides
exactly with the subgraph above. By Theorem 5.3.1, we obtain ${\rm\textbf{Z}}(G)\cap N = 1$ and one of the assertions 1, 3 and 5 follow.

b) Suppose that $N/({\rm \textbf{Z}}(G)\cap N)$ is a CP-group. By Theorem 2.4, its order is divisible by at most two primes. For the rest of the proof, we exclude the case in which $N$ has prime power
order, otherwise, we have case 1. Notice that if ${\rm\textbf{Z}}(G)\cap N = 1$, then by applying Theorems 6.2 and 2.4, we obtain case 5. Therefore, by applying Lemma 4.2.1, we can also assume for the
rest of the proof that $|{\rm\textbf{Z}}(G)\cap N|_2=2$. We distinguish two subcases depending on the connectivity of $\Gamma_G(N)$.

b.1) Suppose that $\Gamma_{G}(N)$ is connected. If it consists only of one vertex, $N$ would be a $p$-group by Theorem 3.1, but we are excluding this case. If $\Gamma_{G}(N)$ has two vertices joined
by an edge, by applying Theorem 4.2.2, we obtain case 2. In any other case, $\Gamma_{G}(N)$ always contains a subgraph of the following way:

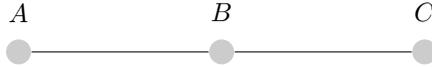
\begin{center}
\begin{tikzpicture}
  [scale=.9,auto=left,every node/.style={circle,fill=black!20}]
 \node [label=$A$](n1) at (0,0) {};
  \node [label=$B$](n2) at (3,0)  {};
  \node [label=$C$](n3) at (6,0)  {};
\foreach \from/\to in {n1/n2,n2/n3}
\draw (\from) -- (\to);
\end{tikzpicture}
\end{center}
 and this turns to a contradiction. Observe that $A$ and $C$ are the unique $G$-classes with cardinalities  $|A|$ and $|C|$, respectively, because otherwise there would exist triangles in
 $\Gamma_G(N)$. If we take a nontrivial $2$-element $z \in {\rm\textbf{Z}}(G)\cap N$, then $zA=A$ and $zC=C$, so  $|\langle z\rangle|$ divides $|A|$ and $|C|$,  which is a contradiction.

b.2) Assume that $\Gamma_{G}(N)$ is disconnected. We distinguish all the possibilities by taking into account Theorem 2.3. First, if $\Gamma_{G}(N)$ has the form:
\begin{center}
\begin{tikzpicture}
  [scale=.9,auto=left,every node/.style={circle,fill=black!20}]
  \node (n1) at (0,0) {};
  \node (n2) at (3,0)  {};
\end{tikzpicture}
\end{center}
 we have  ${\rm\textbf{Z}}(G)\cap N = 1$ by Lemma 4.1.1. Also, by Theorem 4.1.2, $N$ is a $\{p,q\}$-group and $N$ satisfies a particular case of 5. Secondly, if $\Gamma_{G}(N)$ has three $G$-classes,
 the graph is:
\begin{center}
\begin{tikzpicture}
  [scale=.9,auto=left,every node/.style={circle,fill=black!20}]
  \node (n1) at (0,0) {};
  \node (n2) at (3,0)  {};
  \node (n3) at (6,0)  {};
\foreach \from/\to in {n1/n2,n3}
\draw (\from) -- (\to);
\end{tikzpicture}
\end{center}
By Theorem 5.1.1, $N$ is a $\{p,q\}$-group and since $N$ is not a $p$-group, then $N$ is quasi-Frobenius with abelian kernel and complement, so we get case 4.\\

The only remaining case is when the structure of $\Gamma_{G}(N)$ is:
\begin{center}
\begin{tikzpicture}
  [scale=.9,auto=left,every node/.style={circle,fill=black!20}]
  \node (n1) at (0,0) {};
  \node (n2) at (3,0)  {};
  \node (n3) at (6,0)  {};
  \node (n4) at (9,0) {};
\foreach \from/\to in {n1/n2,n3/n4}
\draw (\from) -- (\to);
\end{tikzpicture}
\end{center}

Then by Theorem 2.2, we know that either $N=P\times A$ with $A\neq 1$ a central $p'$-subgroup of $G$, or $N$ is quasi-Frobenius with kernel $K$ and complement $H$, both of them abelian. In the former
case, since the $G$-classes of $N$ are equal to the product of the $G$-classes of $P$ by the classes  of $A$,  necessarily we have $A \cong \mathbb{Z}_2$.

 We assume then the second possibility, that is,  $N$ is quasi-Frobenius with kernel $K$ and complement $H$, both abelian. As $N/({\bf Z}(G)\cap N)$ is a CP-group, it certainly follows that $Z:=
 {\rm\textbf{Z}}(G)\cap N={\rm \textbf{Z}}(N)$. We have proved above that $N/Z$ is a $\{p,q\}$-group and now we prove that $Z \cong \mathbb{Z}_2$. If $z \in Z$ is an $r$-element with $r$ an odd prime,
 then $r$ needs to be distinct from $p$ or $q$, say for instance $p \neq r$. We take a $p$-element $x \in N \setminus Z$ and then the $G$-classes $x^G$, $(xz)^G$ and $(xz^2)^G$ form a triangle in
 $\Gamma_G(N)$, a contradiction. We conclude that $Z \cong \mathbb{Z}_2$, as wanted. Finally, we prove that either $p$ or $q$ has to be $2$, that is, $N$ is a $\{2,p\}$-group (or a $\{2,q\}$-group),
 and we obtain case 4 and the theorem is finished. Suppose that $p$ and $q$ are odd primes.  It easy to see that the graph $\Gamma_{G/Z}(N/Z)$ has exactly two vertices without edge, and by using
 Theorem 4.1.1, we have $|K/Z|=p^n$, $|H/Z|=q$, with $n\geq 1$. We take a $p$-element $a\in K$, a $q$-element $b\in H$ and $1\neq z \in {\rm\textbf{Z}}(G)\cap N$.  It follows that the four vertices of
 $\Gamma_G(N)$ are exactly $a^G$, $(az)^G$, $b^G$ and $(bz)^G$.  Moreover, observe that $a^G \cup (az)^G =K \setminus Z$ and $b^G \cup (bz)^G= N\setminus K$. By taking cardinalities, we have
$$2|a^ G|= 2p^n-2 \quad {\rm and} \quad 2|b^G|= 2p^nq-2p^n,$$
so $|a^G|= p^n-1$ and $|b^G|=p^n (q-1)$. However, both numbers are even because $p$ and $q$ are odd, and this leads to a contradiction.  $\Box$\\

\end{document}